\documentclass[12pt,a4paper]{amsart}

\usepackage{amsfonts, amsmath, amssymb, amsthm, hyperref}

\usepackage{a4wide}

\newtheorem{thm}{Theorem}
\newtheorem{lem}{Lemma}
\newtheorem{con}{Conjecture}

\theoremstyle{definition}
\newtheorem{defn}{Definition}

\theoremstyle{remark}

\newcommand{\pt}{\mathrm{pt}}

\DeclareMathOperator{\Ind}{Ind}

\renewcommand{\int}{\mathop{\rm int}}

\renewcommand{\epsilon}{\varepsilon}

\begin{document}

\title{Knaster's problem for almost $(Z_p)^k$-orbits}

\author{R.N.~Karasev}
\thanks{The research of R.N.~Karasev was supported by the President's of Russian Federation grant MK-1005.2008.1, and partially supported by the Dynasty Foundation.}

\email{r\_n\_karasev@mail.ru}
\address{
Roman Karasev, Dept. of Mathematics, Moscow Institute of Physics
and Technology, Institutskiy per. 9, Dolgoprudny, Russia 141700}

\author{A.Yu.~Volovikov}
\thanks{The research of A.Yu.~Volovikov was partially supported by the Russian Foundation for Basic Research grant No. 08-01-00663}

\email{a\_volov@list.ru}
\address{
Alexey Volovikov, Department of Mathematics, University of Texas at Brownsville, 80
Fort Brown, Brownsville, TX, 78520, USA}

\address{
Alexey Volovikov, Department of Higher Mathematics, Moscow State Institute of
Radio-Engineering, Electronics and Automation (Technical
University), Pr. Vernadskogo 78, Moscow 117454, Russia}


\begin{abstract}
In this paper some new cases of Knaster's problem on continuous maps from spheres are established. In particular, we consider an almost orbit of a $p$-torus $X$ on the sphere, a continuous map $f$ from the sphere to the real line or real plane, and show that $X$ can be rotated so that $f$ becomes constant on $X$.
\end{abstract}

\maketitle

\section{Introduction}

In~\cite{kna1947} the following conjecture (Knaster's problem) was formulated.

\begin{con}
\label{knaster}
Let $S^{d-1}$ be a unit sphere in $\mathbb R^d$. Suppose we are given $m = d-k+1$ points $x_1, \ldots, x_m\in S^{d-1}$ and a continuous map $f: S^{d-1}\to \mathbb R^k$. Then there exists a rotation $\rho\in SO(d)$ such that
$$
f(\rho(x_1)) = f(\rho(x_2)) = \dots = f(\rho(x_m)).
$$
\end{con}

In papers~\cite{kasha2003,hiri2005} it was shown that for certain sets $\{x_1, \ldots, x_m\}\subset S^{d-1}$ Knaster's conjecture fails, such counterexamples exist for every $k>2$, for $k=2$ and $d\ge 5$, for $k=1$ and $d\ge 67$. 

Still it is possible to prove Knaster's conjecture in some particular cases of sets. In~\cite{vol1992} the set of points was some orbit of the action of a $p$-torus $G=(Z_p)^k$ on $\mathbb R[G]$ for $k=1$ and on $\mathbb R[G]\oplus\mathbb R$ for $k=2$. Here we prove some similar results, the set of points being a $(Z_p)^k$-orbit minus one point. 

The group algebra $\mathbb R[G]$ is supposed to have left $G$-action, unless otherwise stated. Considered as a $G$-representation, $\mathbb R[G]$ may have a $G$-invariant inner product. In fact, the space of invariant inner products has the dimension equal to the number of distinct irreducible $G$-representations in $\mathbb R[G]$ (for a commutative $G$), for a $p$-torus $G=(Z_p)^k$ the dimension of this space is $\frac{p^k+1}{2}$ for odd $p$, and $p^k$ for $p=2$.

\begin{defn}
Denote $I[G]\subset\mathbb R[G]$ the $G$-invariant subspace in $\mathbb R[G]$ consisting of
$$
\sum_{g\in G} \alpha_g g,\quad\text{with}\quad \sum_{g\in G} \alpha_g = 0.
$$
Note that its orthogonal complement (w.r.t. any $G$-invariant inner product) is the one-dimensional space with trivial $G$-action.
\end{defn}

In the sequel we consider a $p$-torus $G=(Z_p)^k$ and denote $q=p^k$.

\begin{thm}
\label{knaster1}
Let $S^{q-2}$ be the unit sphere of $I[G]$ w.r.t. some $G$-invariant inner product, denoted by $(\cdot, \cdot)$. Then conjecture~\ref{knaster} holds for $k=1$, the rotations w.r.t. $(\cdot, \cdot)$, and the set $Gx\setminus\{x\}$, where $x\in S^{q-2}$ is any point.
\end{thm}

\begin{thm}
\label{knaster2}
Let $S^{q-1}$ be the unit sphere of $\mathbb R[G]$ w.r.t. some $G$-invariant inner product $(\cdot, \cdot)$. Then conjecture~\ref{knaster} holds for $k=2$, $q$ odd, the rotations w.r.t. $(\cdot, \cdot)$, and the set $Gx\setminus\{x\}$, where $x\in S^{q-1}$ is any point.
\end{thm}

In fact, the last theorem may be formulated a little stronger. For example, Theorem~\ref{knaster2d} (see below) gives the following statement. Let $x\in S^{q-1}$ be as in the theorem, and let $f_1, f_2 : S^{q-1}\to \mathbb R$ be two continuous functions. Then for some rotation $\rho$ and two constants $c_1, c_2$
\begin{eqnarray*}
&\forall g\in G\ &f_1(\rho(gx)) = c_1\\
&\forall g\in G,\ g\not=e\ &f_2(\rho(gx))=c_2.
\end{eqnarray*}

\section{Equivariant cohomology of $G$-spaces}
\label{eq-cohomology}

We consider topological spaces with continuous action of a finite group $G$ and continuous maps between such spaces that commute with the action of $G$. We call them $G$-spaces and $G$-maps.

Let us consider the group $G=(Z_p)^k$ and list the results (mostly from~\cite{vol2005}) that we need in this paper.

The cohomology is taken with coefficients in $Z_p$, in the notation we omit the coefficients. 

Consider the algebra of $G$-equivariant (in the sense of Borel) cohomology of the point $A_G = H_G^*(\pt) = H^*(BG)$. For any $G$-space $X$ the natural map $X\to\pt$ induces the natural map of cohomology $\pi_X^* : A_G\to H_G^*(X)$.

For a group $G=(Z_p)^k$ the algebra $A_G$ (see~\cite{hsiang1975}) has the following structure. For odd $p$, it has $2k$ multiplicative generators $v_i,u_i$ with dimensions $\dim v_i = 1$ and $\dim u_i = 2$ and relations
$$
v_i^2 = 0,\quad\beta{v_i} = u_i.
$$
Here we denote $\beta(x)$ the Bockstein homomorphism.

For a group $G=(Z_2)^k$ the algebra $A_G$ is the algebra of polynomials of $k$ one-dimensional generators $v_i$.

The powerful tool of studying $G$-spaces is the following spectral sequence (see~\cite{hsiang1975,mcc2001}).

\begin{thm}
\label{specseqeq}
There exists a spectral sequence with $E_2$-term
$$
E_2^{x, y} = H^x(BG, \mathcal H^y(X, Z_p)),
$$
that converges to the graded module, associated with the filtration of $H_G^*(X, Z_p)$.

The system of coefficients $\mathcal H^y(X, Z_p)$ is obtained from the cohomology $H^y(X, Z_p)$ by the action of $G = \pi_1(BG)$. The differentials of this spectral sequence are homomorphisms of $H^*(BG, Z_p)$-modules.
\end{thm}

For every term $E_r(X)$ of this spectral sequence there is a natural map $\pi^*_r : A_G\to E_r(X)$. 

\begin{defn}
Denote the kernel of the map $\pi^*_r$ by $\Ind^r_G X$. 
\end{defn}

The ideal-valued index of a $G$-space was introduced in~\cite{fh1988}, the above filtered version was introduced in~\cite{vol2000}. Remind the properties of $\Ind^r_G X$, that are obvious by the definition. We omit the subscript $G$ when a single group is considered.

\begin{itemize}
\item 
If there is a $G$-map $f:X\to Y$, then $\Ind^r X\supseteq \Ind^r Y$.

\item
$\Ind^{r+1} X$ may differ from $\Ind^r X$ only in dimensions $\ge r$.

\item
$\bigcup_r \Ind^r X = \Ind X = \ker \pi_X^* : A_G\to H_G^*(X)$.
\end{itemize}

The first property in this list is very useful to prove nonexistence of $G$-maps. Following~\cite{vol2005} we define a numeric invariant of this ideal filtering $\Ind_G^r X$.

\begin{defn}
Put 
$$
i_G(X) = \max \{r : \Ind_G^r X = 0\}.
$$
\end{defn}

It is easy to see that $i_G(X)\ge 1$ for any $G$-space $X$, $i_G(X)\ge 2$ for a connected $G$-space $X$, and $i_G(X)$ may be equal to $+\infty$. Moreover, for a $G$-space $X$ without fixed points, $G$-homotopy equivalent to a finite $G$-$CW$-complex, $i_G(X)\le \dim X + 1$.

From the definition of $\Ind_G^r X$ it follows that if there exists a $G$-map $f:X\to Y$, then $i_G(X)\le i_G(Y)$ (the monotonicity property). 


The definition of $i_G(X)$ can be further extended.

\begin{defn}
Define the index of a cohomology class $\alpha\in A_G$ on a $G$-space $X$
$$
i_G(\alpha, X) = \max\{r: \alpha\not\in \Ind_G^r X\}.
$$
It may equal $+\infty$ if $\alpha\not\in \Ind_G X$
\end{defn}

It is clear from the definition that either $i_G(\alpha, X)=+\infty$, or $i_G(\alpha, X)\le \dim \alpha$ and $i_G(\alpha, X)\le \dim X + 1$ (for a finite $G$-$CW$-complex). Moreover, for any $G$-map $f : X\to Y$ we have the monotonicity property
$$
i_G(\alpha, X)\le i_G(\alpha, Y).
$$

\section{Reformulations}
\label{reformulations}

We reformulate Theorems~\ref{knaster1} and \ref{knaster2} in a more general way.

\begin{thm}
\label{knaster1d}
Let $S^{q-2}$ be the unit sphere of $I[G]$ w.r.t. some $G$-invariant inner product, and let $f:S^{q-2}\to\mathbb R$ be some continuous function. Consider $x\in S^{q-2}$, the vector $v=\sum_{g\in G} g\in\mathbb R[G]$ and some other vector $w\in\mathbb R[G]$, non-collinear to $v$.

Then for some rotation $\rho\in SO(q-1)$ the vector $\sum_{g\in G}f(\rho(gx))g\in\mathbb R[G]$ is in the linear span of $v$ and $w$.
\end{thm}

Theorem~\ref{knaster1} follows from this theorem in the following way. Put $w=e\in\mathbb R[G]$. Then by Theorem~\ref{knaster1d} there exists a rotation $\rho$ such that for some $\alpha,\beta\in\mathbb R$
$$
\forall g\in G,\ g\not=e,\ f(\rho(gx)) = \alpha,\quad f(\rho(x)) = \alpha +\beta.
$$
That is exactly the statement of Theorem~\ref{knaster1}.

\begin{thm}
\label{knaster2d}
Let $S^{q-1}$ be the unit sphere of $\mathbb R[G]$ w.r.t. some $G$-invariant inner product, and let $f:S^{q-1}\to\mathbb R^2$ be some continuous map with coordinates $f_1, f_2$. Let $q$ be odd. Consider $x\in S^{q-1}$, the vectors $v=\sum_{g\in G}g\oplus 0\in\mathbb R[G]\oplus\mathbb R[G]$, $u=0\oplus\sum_{g\in G} g\in\mathbb R[G]\oplus\mathbb R[G]$ and some other vector $w\in\mathbb R[G]\oplus\mathbb R[G]$, non-coplanar to $v, u$.

Then for some rotation $\rho\in SO(q)$ the vector 
$$
\sum_{g\in G}f_1(\rho(gx))g\oplus\sum_{g\in G}f_2(\rho(gx))g\in\mathbb R[G]\oplus\mathbb R[G]
$$ 
is in the linear span of $v, u, w$.
\end{thm}

Again, Theorem~\ref{knaster2} (and its stronger version in the remark after Theorem~\ref{knaster2}) follows from this theorem by taking a vector $w=e\oplus 0$, similar to the previous remark.

\section{Proof of Theorem~\ref{knaster1d} in the case of odd $q$}

\label{knaster1d-odd}

In this section $q=p^k$, $p$ is an odd prime, $G=(Z_p)^k$. Define for any $\rho\in SO(q-1)$
$$
\phi(\rho) = \sum_{g\in G} f(\rho(g(x)))g\in\mathbb R[G].
$$
For any $h\in G$ we have
$$
\phi(\rho\circ h^{-1}) = \sum_{g\in G} f(\rho(h^{-1}g(x)))g = \sum_{g\in G} f(\rho(h^{-1}g(x)))hh^{-1}g=\sum_{g\in G} f(\rho(g(x)))hg.
$$
Thus the map $\phi : SO(q-1)\to \mathbb R[G]$ is a $G$-map for the left action of $G$ on $SO(q-1)$ by right multiplications by $g^{-1}\in G$, and for the standard left action of $G$ on $\mathbb R[G]$. 

Denote for any $g\in G$ by $L_g=(v, gw)\subset\mathbb R[G]$ the $2$-dimensional subspaces. Assume the contrary: that is the image of $\phi$ does not intersect $\bigcup_{g\in G} L_g$. So $\phi$ maps $SO(q-1)$ to the space $Y=\mathbb R[G]\setminus \bigcup_{g\in G} L_g$. The natural projection $\pi : Y\to\mathbb R[G]/(v)=V$ gives a homotopy equivalence between $Y$ and $V\setminus \bigcup_{g\in G} \mathbb R\pi(gw)$, the latter space is homotopically $q-2$-dimensional sphere without several points, hence it is a wedge of $q-3$-dimensional spheres. $G$ acts on $Y$ without fixed points, so $i_G(Y)\le q-2$. 

In~\cite{vol1992} it was shown that $i_G(SO(q-1)) = q-1$ w.r.t. the considered $G$-action. Here we give a short explanation. In the spectral sequence of Theorem~\ref{specseqeq} all multiplicative generators of $H^*(SO(q-1), Z_p)$ are transgressive, because they are pullbacks of the transgressive generators of $H^*(SO(q-1), Z_p)$ in the spectral sequence of the fiber bundle $\pi_{SO(q-1)} : ESO(q-1)\to BSO(q-1)$. So the first nonzero $\Ind_G^r SO(q-1)$ corresponds to the first nonzero characteristic class of the $G$-representation $I[G]$ in the cohomology ring $A_G$. It was shown in~\cite{vol1992} that this is the Euler class of $I[G]$ of dimension $q-1$.

So we have a contradiction with the monotonicity of $i_G(X)$.

\section{Proof of Theorem~\ref{knaster2d}}
\label{knaster2d-proof}

Similar to the previous proof, we consider the $G$-map $\phi: SO(q)\to \mathbb R[G]\oplus\mathbb R[G]$, given by the formula
$$
\phi(\rho) = \sum_{g\in G} f_1(\rho(g(x)))g\oplus \sum_{g\in G} f_2(\rho(g(x)))g\in\mathbb R[G]\oplus\mathbb R[G].
$$ 

Take the composition $\psi = \pi\cdot \phi$ with the projection $\pi:\mathbb R[G]\oplus\mathbb R[G]\to I[G]\oplus I[G] = V$. Assume the contrary: that is the map $\phi$ does not intersect the linear span of $u$ and $v$ in $\mathbb R[G]\oplus\mathbb R[G]$ and $\psi$ does not intersect the linear span of $gw$ for any $g\in G$ in $V$, it means that the image of $\psi$ is in the space $Y = V\setminus \bigcup_{g\in G} \mathbb R\pi(gw)$. 

Let $e\in A_G$ be the Euler class of $V$. From the spectral sequence of Theorem~\ref{specseqeq} it is obvious that $i_G(e, V\setminus\{0\}) = 2q-2$, because the spectral sequence for $V\setminus\{0\}$ has the only nontrivial differential that kills the Euler class $e$. Since $Y\subset V\setminus\{0\}$, then $i_G(e, Y) < +\infty$. Similar to the previous proof, the space $Y$ is homotopically a wedge of $2q-4$-dimensional spheres, so $i_G(e, Y) \le \dim Y + 1 = 2q-3$. 

In~\cite{vol1992} it was shown that $i_G(e, SO(q)) = 2q-2$, because $e$ is in the image of the transgression in the spectral sequence and $e$ is not contained in the ideal of $A_G$, generated by the characteristic classes of $SO(q)$ of lesser dimension. So we again have a contradiction with the monotonicity of $i_G(e, X)$.

\section{Proof of Theorem~\ref{knaster1d} in the case of even $q$}

In this section $q=2^k$, $G=(Z_2)^k$. We use the notation from the odd case in Section~\ref{knaster1d-odd}. Note that the case $q=2$ is trivial, and if $q\ge 4$ then $G$ acts on $I[G]$ by transforms with positive determinant, so the group $SO(q-1)$ can be considered as the configuration space.

Assume the contrary: the image $\phi(SO(q-1))$ is in $Y=\mathbb R[G]\setminus \bigcup_{g\in G} L_g$.

Denote the Stiefel-Whitney classes of $I[G]$ in $A_G$ by $w_k$. We need the following lemma, stated in~\cite{vol1992}, based on results from~\cite{dic1911,mui1975}.

\begin{lem}
\label{euler-nz}
The only nonzero Stiefel-Whitney classes of $I[G]$ are $w_{q-2^l}\in A_G$ $(l=0,\ldots, k)$, the classes $w_{q-2^l}$ $(l=0,\ldots, k-1)$ are algebraically independent and form a regular sequence, hence $w_{q-1}$ is nonzero and not contained in the ideal of $A_G$, generated by $w_k$ with $k<q-1$. 
\end{lem}

Similar to the proof of Theorem~\ref{knaster2d} in Section~\ref{knaster2d-proof}, we find that $i_G(w_{q-1}, Y)\le \dim Y + 1 = q-2$.

Now we apply the spectral sequence of Theorem~\ref{specseqeq} to the $G$-space $SO(q-1)$. The action of $G$ on $SO(q-1)$ is the restriction of action of $SO(q-1)$ on itself, the latter group being connected, hence $G$ acts trivially on $H^*(SO(q-1), Z_2)$.

The results of~\cite{bor1953} imply that the differentials in this spectral sequence are generated by transgressions that send the primitive (in terms of~\cite{bor1953}) elements of $H^*(SO(q-1), Z_2)$ to the Stiefel-Whitney classes $w_k$ (see Proposition~23.1 in~\cite{bor1953}). Thus Lemma~\ref{euler-nz} implies that $i_G(w_{q-1}, SO(q-1)) = q-1$, and the existence of the $G$-map $\phi$ contradicts the monotonicity of $i_G(w_{q-1}, X)$.

\end{document}